\documentclass[12pt]{article}
\usepackage{amsmath,amsfonts,amssymb,amsthm}
\newcommand{\1}{{\bf 1}}
\newcommand{\id}{{\rm id}}
\newcommand{\Hom}{{\rm Hom}\,}
\newcommand{\End}{{\rm End}\,}
\newcommand{\Res}{{\rm Res}\,}

\newcommand{\gr}{{\rm gr}}
\newcommand{\QP}{{\mathcal{Q}}}
\newcommand{\Prim}{{\mathcal{P}}}
\newcommand{\har}{{\rm span}}
\newcommand{\divi}{\operatorname{ord}}
\newcommand{\F}{{\mathcal F}}
\newcommand{\V}{{\mathcal V}}
\newcommand{\qV}{{\mathcal{QV}}}

\newcommand{\g}{\mathfrak{g}}
\newcommand{\Vir}{\mathcal{V}ir}

\newcommand{\GG}{{\mathcal{Q}}}
\newcommand{\PG}{{\mathcal{G}}}
\newcommand{\Z}{\mathbb{Z}}

\newcommand{\C}{\mathbb{C}}
\newcommand{\N}{\mathbb{N}}

\newcommand{\h}{\mathfrak{h}}

\newcommand{\TS}{\widetilde{\Sigma}}
\newcommand{\wt}[1]{|#1|}

\newtheorem{theorem}{Theorem}[section]
\newtheorem{proposition}[theorem]{Proposition}
\newtheorem{lemma}[theorem]{Lemma}

\theoremstyle{definition}
\newtheorem{definition}[theorem]{Definition}
\newtheorem{example}[theorem]{Example}
\newtheorem{note}[theorem]{Note}

\theoremstyle{remark}

\numberwithin{equation}{section}

\catcode`\@=11
\renewcommand\section{\@startsection{section}{1}{\z@}%
                                      {-3.25ex\@plus -1ex \@minus -.2ex}%
                                      {1.5ex \@plus .2ex}%
                                      {\normalfont\large\bfseries}}
\catcode`\@=\active
\catcode`\@=11
\renewcommand\subsection{\@startsection{subsection}{2}{\z@}%
                                      {-3.25ex\@plus -1ex \@minus -.2ex}%
                                      {1.5ex \@plus .2ex}%
                                      {\normalfont\large\bfseries}}
\catcode`\@=\active

\begin{document}

\begin{large}
\begin{center}
\textbf{Finiteness of conformal blocks over compact Riemann surfaces}
\end{center}
\end{large}

\vskip 2ex
\begin{center}
Toshiyuki Abe\footnote{Supported by JSPS Research Fellowships for Young 
Scientists.} and
Kiyokazu Nagatomo\footnote{Supported in part by  Grant-in-Aid for 
Scientific Research, Japan Society for the Promotion of Science.} \\

\vskip 2ex
Department of Mathematics, Graduate School of Science\\
  Osaka University, Toyonaka, Osaka 560-0043, Japan
\end{center}

\vskip 3ex
\noindent
\begin{small}
\textbf{Abstract:}
We study conformal blocks (the space of correlation functions) over compact Riemann surfaces
associated to vertex operator algebras which are the sum of highest weight modules for the
underlying Virasoro algebra. Under the fairly  general condition, for instance, $C_2$-finiteness, 
we prove that conformal blocks are of finite dimensional. This, in particular, shows the finiteness
of conformal blocks  for many well-known conformal field theories including WZNW model
and the minimal model.
\end{small}

\baselineskip 3ex
\section*{Introduction}
In \cite{AN} we have shown that conformal blocks over the \textit{projective
line} associated to a vertex operator algebra(VOA) $V$ are of finite dimensional if modules for
$V$ satisfy some finiteness condition. In this paper we generalize these results
to conformal blocks over any \textit{compact Riemann surfaces}. 
More precisely we will prove that if $V$-modules of our concern as well as $V$ are
$C_2$-finite then corresponding conformal blocks are of finite dimensional.  
The main reason that we need $C_2$-finiteness of $V$ in this case is caused by Weierstrass gaps, i.e., 
we are not able to find  meromorphic differentials with poles of some exceptional orders.

%
%

Though in this paper the notion of conformal blocks are defined in a purely mathematical way, 
the definition goes bock to the notion of correlation functions in conformal field theory (CFT) initiated
by \cite{BPZ}. CFT's are supposed to have at least two properties, one of which is the finiteness of
conformal blocks, and the other is the factorization property; the latter enables us to determine the dimension of
conformal blocks by  fusion rules (the space of
$3$-point correlation functions or its dimension).  Like as other objects in physics every CFT has it own
symmetry group (Lie algebra): affine Lie algebras for WZNW model and the Virasoro algebra for the minimal
model, for instance. We will study  ``general'' CFT's, where ``general'' means that
the symmetry is described by a VOA. Such CFT's were first proposed and studied by Zhu
\cite{Z2}, however two main issues, i.e., finiteness of conformal blocks and the factorization theorem  of these CFT's were left open. 

We should point out two main differences between our general CFT's and  known CFT's.  Conformal blocks are the space of
correlation functions of primary fields and the Virasoro  fields. On the one hand, for instance, in WZNW
model,  the space of states is generated by currents which are all conformal weight $1$
primary fields, and the Virasoro field of the theory, which has conformal weight $2$, is given by Sugawara construction
using the  currents, and so we only need information on meromorphic functions on the Riemann surface 
to study conformal blocks. However, in general, we  have primary fields of higher conformal weight
$\Delta$, and we have to know the geometry of the line bundle $\kappa^{1-\Delta}$ where $\kappa$ is  the
canonical line bundle. The minimal model is generated by a conformal weight $2$ field (the Virasoro field)
and the analysis of the line bundle 
$\kappa^{-1}$ is necessary, though it is still not so complicated.

A part of ideas in \cite{TUY} was generalized to VOA's in \cite{Z2}: Zhu first generalized the notion of currents and the energy-momentum
field to the notion of \textit{global vertex operators} associated to any primary states and then gave a very general definition of conformal
blocks; more precisely, conformal blocks are defined in the almost same way as \cite{TUY}. However this definition using primary fields and the
Virasoro field is not convenient because Fourier modes of primary fields and the Virasoro field do not form a Lie algebra. Zhu then introduced
so-called \textit{quasi-global vertex operators} which are defined by using quasi-primary states. 
The quasi-global vertex operators now form a Lie algebra under the fairly general assumption for $V$. 
The point is that we can characterize conformal blocks in terms of quasi-global vertex operators, 
which is one of the main results in \cite{Z2}.

In many examples of CFT's the key fact for the finiteness of conformal blocks is the finite dimensionality of the space of coinvariants; several
examples  are known such as WZNW model and the minimal model.  The notion in VOA theory corresponding to the finiteness of ``coinvariants''  is
$C_2$-finiteness condition introduced in \cite{Z1}. Using the notion of Frenkel-Zhu bimodules \cite{FZ} the finiteness of fusion rules is proved
in \cite{Li} for $C_2$-finite modules; more precisely, the weaker condition called $B_{1}$-finite is enough for the finiteness of fusion rules. 

In this paper we prove the finiteness of conformal blocks over pointed compact Riemann surfaces 
associated to $C_{2}$-finite vertex operator algebras which are sum of highest weight modules of the Virasoro algebra and $B_{1}$-finite
$V$-modules; we should mention that the our notion of $B_1$-finiteness is different from Li's one.
The method used here basically follows \cite{TUY}, while we work on fairly general set up.  
The proof of finiteness of conformal blocks over compact Riemann surfaces is reduced to find 
a non-trivial meromorphic sections with poles of specified positions and orders. The main difference between the case of the projective line
\cite{AN} and the case studied in the paper is that on a given Riemann surface we are not always able to find a meromorphic form which has poles
at the prescribed points and orders; this point will be elaborated in the paper.

The conformal blocks over pointed projective lines are also studied in \cite{GN}. 
The definition of conformal blocks looks different from the one of \cite{Z2} and \cite{AN},
but their method has a great influence to our work.

In section 1 we review some basics for vertex operator algebras.
In section 2 we recall the notion of $C_{n}$-finiteness condition $(n\geq2)$ and $B_{1}$-finiteness condition for modules,
and state several known results  concerning these finiteness conditions. 
The notion of conformal blocks is defined in section 3 following \cite{Z2}; conformal blocks are built of a
compact Riemann surface $\Sigma$, a finite set $A$ of distinct points on $\Sigma$  and
$V$-modules
$W^i\,(i\in A)$.  We introduce a filtration on the Lie algebra $\GG(\TS)$ of all quasi-global
vertex operators on a compact Riemann surface
$\Sigma$. The Lie algebra
$\GG(\TS)$ acts on the tensor product vector space $W_A=\bigotimes_{i\in A}W^i$, on which we define a filtration 
so that $W_A$ becomes a filtered $\GG(\TS)$-module. These filtrations have their origins in
\cite{TUY}, but of course they are appropriately generalized to fit our purposes.  Section 4 is the core of this paper, in which
we prove that conformal blocks associated to a $C_{2}$-finite quasi-primary generated vertex operator algebra
and $B_1$-finite $V$-modules are of finite dimensional; we prove the main theorem 
by showing the existence of a surjective map $\bigotimes_{i\in A}W^i/B_1(W^i)\rightarrow \gr_\bullet W_A/\GG(\TS)W_A$. 
 In particular if all modules $W^i\,(i\in A)$ are 
$C_2$-finite then the corresponding conformal block is of finite dimensional. 
Finally in section 5 we discuss several examples of $C_2$-finite vertex operator algebras whose irreducible modules are $B_1$-finite.

After we finished the work we learned Buhl's result \cite{Bu} that any finitely generated weak module for a $C_2$-finite (which is called $C_2$
co-finite in \cite{Bu}) vertex operator algebra is $C_n$-finite for all $n\geq2$. Therefore our $B_1$-finiteness assumption for modules is
not necessary.

\section{Vertex operator algebras and their modules}
Let $V$ be a vertex operator algebra with the vacuum element $\1$ and the Virasoro element $\omega$
(see \cite{FLM}, \cite{FHL}, \cite{MN}), i.e., the vector space $V$ is equipped with countably many bilinear operation
$(a,b)\mapsto a(n)b\,(a,\,b\in V)$ for any integer $n$.
For any $a\in V$, we denote the vertex operator associated to $a$ by $Y(a,x)=\sum_{n\in\Z}a(n)x^{-n-1}$ where $a(n)\in \End(V)$ is defined by
$b\mapsto a(n)b$ for all $b\in V$. The operators $L_n:=\omega(n+1)\,(n\in\Z)$ form a representation of the Virasoro algebra on $V$,
and the vector space $V$ is
$\N$-graded with the grading operator $L_0$, i.e., $V=\bigoplus_{n=0}^\infty V(n),\, L_0|V(n)=n\,\id$.
The operator $L_{-1}$ is assumed to satisfy $\frac{\partial}{\partial x}Y(a,x)=Y(L_{-1}a,x)$, i.e., $-na(n-1)=(L_{-1}a)(n)$ for all $a\in V$ and
$n\in\Z$.

An element $a\in V$ satisfying $L_{1}a=L_{2}a=0$ is called a {\it primary vector}, while an element satisfying only $L_1a=0$ is called a
\textit{quasi-primary vector}.  Let $\Prim(V)$ and $\QP(V)$ be the sets of all primary and quasi-primary vectors, respectively.
We see that those two vector subspaces of $V$ are graded, i.e., 
$\Prim(V)=\bigoplus_{n=0}^{\infty}\Prim(V)\cap V(n)$ and $\QP(V)=\bigoplus_{n=0}^{\infty}\QP(V)\cap V(n)$.

\begin{definition}\label{Def1}
A vertex operator algebra $V$ satisfying $V=\sum_{k=0}^\infty L_{-1}^{k} \QP(V)$ is called \textit{quasi-primary generated}.
\end{definition}
It is known that $V$ is quasi-primary generated if and only if $V_{1}\subset \QP(V)$ (\cite{DLinM}).
For $V$ quasi-primary generated we see that $V=\bigoplus_{k=0}^\infty L_{-1}^k\QP(V)$ if and only if $L_{-1}V(0)=0$.

\begin{definition}\label{Def2}{\rm
A {\em weak $V$-module} is a vector space $W$ equipped with a linear map
\begin{align*}
Y_{W}:V&\to (\End W)[[x,x^{-1}]]\\
a&\mapsto Y_{W}(a,x)=\sum_{n\in\Z} a(n)x^{-n-1},\ (a(n)\in \End W)
\end{align*}
which satisfies the following conditions for all $a,\,b\in V$ and $w\in W$; $Y_{W}(a,x)w\in W((x))$, $Y_{W}(\1,x)=\id_{W}$,
and for all integers $p,\,q,\,r\in\Z$,
\begin{multline}\label{eqn:b1}
\sum_{i=0}^\infty
\binom{p}{i}(a(r+i)b){(p+q-i)}w \\
=\sum_{i=0}^\infty (-1)^i\binom{r}{i}\left(a{(p+r-i)}(b{(q + i)}w) -(-1)^rb{(q+r-i)}(a{(p+i)}w)\right).
\end{multline}
}
\end{definition}

The identity (\ref{eqn:b1}) is equivalent to the set of the following two formulas for $a,\,b\in V$ and $w\in W$ (cf. \cite[\S\,4.3]{MN});
one is called the \textit{associativity formula} 
\begin{multline}\label{asso1}
(a(-n)b)(-q)w\\=\sum_{i=0}^{\infty}\binom{-n}{i}(-1)^{i}(a(-n-i)b(-q+i)w-(-1 
)^{n}b(-n-q-i)a(i)w),
\end{multline}
and the other is called the \textit{commutator formula}
\[
[a(p),b(q)]w=\sum_{i=0}^{\infty}\binom{p}{i}(a(i)b)(p+q-i)w.
\]
Using the commutator formula we see that $L_n\,(n\in\Z)$ form a representation on $W$ for the Virasoro algebra.
By definition of vertex operator algebras $L_{-1}a=(L_{-1}a)(-1)\1=a(-2)\1$ for any  $a\in V$. 
Then the associativity formula for $(a(-2)\1)(q)w$ shows that $(L_{-1}a)(q)w=-qa(q-1)w$ for all $a\in V$ and $w\in W$.

\begin{definition}\label{Def4}{\rm  
A {\em $V$-module} $W$ is a weak $V$-module on which $L_{0}$ acts semisimply, i.e.,
$W=\bigoplus_{\lambda\in\C}W(\lambda),\, L_0|W(\lambda)=\lambda\,\id$, and for fixed
$\lambda\in\C$, $W(\lambda+n)=0$ for all sufficiently large integers $n$. For $\lambda\in\C$, 
a nonzero vector $w$ in $W(\lambda)$ is called {\em homogeneous vector of weight $\lambda$}, and its weight is denoted by $\wt{w}$.
}
\end{definition}
\noindent
Whenever we write $\wt{w}$ the element $w$ is supposed to be homogeneous of weight $\wt{w}$.

\section{Finiteness conditions of vertex operator algebras}
We recall the notion of $C_{n}$-finiteness (see \cite{Z1} for $n=2$, and \cite{Li} for general $n(\geq2)$). 
We review the notion of $B_{1}$-finiteness in \cite{AN}, and state several results;
the most importance is that for a $C_2$-finite vertex operator algebra
$V$ any $B_{1}$-finite weak $V$-module is $C_n$-finite for all $n\geq2$, which was proved in \cite{AN}.
%
%

\begin{definition}\label{Defn3}{\rm (\cite{Z1}, \cite{Li})
For any positive integer $n(\geq2)$ we denote by $C_{n}(W)$ the subspace of $W$,
which is linearly spanned by $a(-n)w$ for all $a\in V$ and $w\in W$. A weak $V$-module $W$ is called {\em $C_{n}$-finite} ($n\geq2$) 
if the vector space $W/C_{n}(W)$ is of finite dimensional.}
\end{definition}

Using $(L_{-1}a)(q)w=-qa(q-1)w$ for all $ a\in V$ and $w\in W$ we see that $C_2(W)\supset C_3(W)\supset \ldots\supset C_n(W)\supset\ldots$,
and that $C_n$-finite module for some $n\geq2$ is $C_2$-finite.
We now let $V$ be a $C_2$-finite vertex operator algebra:


\begin{proposition}\label{Pro7}{\rm (\cite[Proposition 8]{GN})} 
Let $V=\bigoplus_{n=0}^\infty V(n)$ be a vertex operator algebra with $V(0)=\C\1$,
and $U$ be any graded subspace such that $V=U\oplus C_2(V)$. Then $V$ is linearly spanned by the vectors 
\begin{align}\label{equ3}
\alpha_{1}(-n_{1})\alpha_{2}(-n_{2})\cdots\alpha_{k}(-n_{k})\1\mbox{ for all }\alpha_i\in U \mbox{ and } n_{1}>n_{2}>\ldots>n_{k}>0.
\end{align} 
\end{proposition}
%

Let $U$ be a graded subspace such that $V=U\oplus C_2(V)$. 
By Proposition \ref{Pro7} we see that the vectors  (\ref{equ3}) for $k\geq n$ belong to $C_{n}(V)$. 
Suppose that $V$ is $C_2$-finite. Then $U$ is of finite dimensional, and we have:
\begin{proposition}\label{Pro9}{\rm (\cite[Theorem 11]{GN})}
Let $V=\bigoplus_{n=0}^\infty V(n),\,V(0)=\C\1$ be a $C_{2}$-finite vertex operator algebra. 
Then $V$ is $C_{n}$-finite for all $n\geq2$.
\end{proposition}

Let $W$ be a weak $V$-module. We denote by $B_{1}(W)$ the subspace of $W$, which is linearly spanned by 
$a(-1)w$ for all homogeneous $a\in V$ with positive weight, i.e., $|a|>0$ and all $w\in W$; we note that $B_1(W)\supset C_2(W)$.
\begin{definition}\label{Def7}{\rm 
A weak $V$-module $W$ is called \textit{$B_{1}$-finite} if the vector space $W/B_{1}(W)$ is of finite dimensional.}
\end{definition}

\begin{note}
(1) A $B_{1}$-finite weak module is called \textit{quasirational} in \cite{GN}. 

\vskip 1ex
\noindent
(2) The notion of $B_1$-finiteness is slightly different from Li's one \cite{Li}.
\end{note}

By definition $C_n$-finite module for $n\geq 2$ is $B_1$-finite. Conversely we have:
\begin{theorem}[\cite{AN}]\label{Thm3}
Let $V=\bigoplus_{n=0}^\infty V(n),\,V(0)=\C\1$ be a $C_{2}$-finite vertex operator algebra. Then any $B_{1}$-finite weak $V$-module is 
$C_{n}$-finite for all $n\geq2$.  
\end{theorem}
%

\section{Conformal blocks}

We recall the definition of conformal blocks and review some properties of them.  Most of materials in this section except a filtration on
quasi-global vertex operators are taken from
\cite{Z2}.

Let $\Sigma$ be a compact Riemann surface and $\kappa$ be the canonical line bundle on $\Sigma$, 
and let us fix $N$ distinct points  $Q_{1},\,Q_{2},\ldots,Q_{N}$ on $\Sigma$. 
For any integer $n$ we denote  by $\Gamma(\Sigma\,;Q_{1},Q_{2},\ldots,Q_{N};\kappa^{n})$ 
the vector space of global meromorphic sections of $\kappa^n$, which possibly have poles only at $Q_{1},\,Q_{2},\ldots,Q_{N}$.

\begin{definition}\label{Def6}{\rm 
Let $\Sigma$ be a compact Riemann surface and $Q_{1},\,Q_{2},\ldots,Q_{N}$ be distinct points on $\Sigma$. 
Let  $z_{i}$ be local coordinates around the $Q_{i}$ such that $z_{i}(Q_{i})=0$. 
A collection of datum $\bar{\Sigma}=(\Sigma\,;Q_{1},\,Q_{2},\ldots,Q_{N};z_{1},\,z_{2},\ldots,z_{N})$
is called an {\em $N$-pointed Riemann surface}. An $N$-pointed Riemann surface $\bar{\Sigma}$
with a set of $V$-modules $W^{i}$ being attached to each point $Q_{i}$
\begin{align}\label{TS1}
\TS=(\Sigma\,;Q_{1},\,Q_{2},\ldots,Q_{N};z_{1},\,z_{2},\ldots,z_{N};W^{1},\, W^{2},\ldots,W^{N})
\end{align}
is called an {\em $N$-labeled Riemann surface}\,.}
\end{definition}

A system of local coordinates $\{(U_{\alpha},z_{\alpha})\}$ of $\Sigma$ is called a {\em projective structure} on $\Sigma$ 
if transition functions $z_{\beta}\circ z_{\alpha}^{-1}$ are  M{\"o}bius transformations for all $\alpha,\,\beta$ such that $U_\alpha\cap U_\beta\ne\emptyset$; any Riemann surface has a projective structure.
Let $\{(U_{\alpha},z_{\alpha})\}$ be a projective structure of $\Sigma$ and $Q_{1},\,Q_{2},\ldots,Q_{N}$ be distinct points of $\Sigma$.
For each $Q_{i}$ we choose a local coordinate $(U_{\alpha},z_{\alpha})$ 
such that $Q_{i}\in U_{\alpha}$, and define a new coordinate near $Q_i$ by $z_{i}=z_{\alpha}-z_{\alpha}(Q_{i})$. 
Then we obtain an $N$-pointed Riemann surface $\bar{\Sigma}$: such a $\bar{\Sigma}$ is called {\it projective}\;. 
The notion of projective $N$-labeled Riemann surface $\TS$ is defined in a same way.

Let $\bar{\Sigma}=(\Sigma\,;Q_{1},\,Q_{2},\ldots,Q_{N};z_{1},\,z_{2},\ldots,z_{N})$ be an $N$-pointed Riemann surface. 
We will define a Lie algebra $\g(V)_{\bar{\Sigma}}^{out}$ associated to $\bar{\Sigma}$. 
Let $V$ be a vertex operator algebra. We set $\widehat{V}=V\otimes\C((t))$ where $\C((t))$ is the ring of formal Laurent power series. 
It is well known that the commutative associative algebra $\C((t))$ with the derivation $d/dt$ naturally becomes a vertex algebra by 
\[
Y\left(f(t),x\right)g(t)= (e^{x\frac{d}{dt}}f(t))g(t).
\]  
The tensor product $\widehat{V}=V\otimes\C((t))$ has a structure of vertex algebra, which is given by
\[
Y\left(a\otimes f(t),x\right)b\otimes g(t)=Y(a,x)b\otimes (e^{x\frac{d}{dt}}f(t))g(t)
\] 
for all $a\otimes f(t)\,(a\in V,\,f(t)\in\C((t))\,)$.
The translation operator is $D=L_{-1}\otimes\id+\id\otimes\frac{d}{dt}$. We set $\g(V)=\widehat{V}/D\widehat{V}$. Then it is well
known that the $0$-th product on $\widehat{V}$ induces a Lie algebra structure on $\g(V)$.  The important point is that any weak $V$-module
becomes a $\g(V)$-module by $(a\otimes f(t))u=\underset{t=0}{\Res}Y(a,t)f(t)u$.

Let $A=\{1,\,2,\,\ldots,N\}$, and set $\g(V)_{A}=\bigoplus_{i\in A}\g(V)_{(i)}$ where $\g(V)_{(i)}=\g(V)$ is a copy of $\g(V)$. 
We now define a linear map 
\[
j_{\bar{\Sigma}}:\bigoplus_{d=0}^{\infty} V(d\,)\otimes\Gamma(\Sigma\,;Q_{1},\,Q_{2},\ldots,Q_{N};\kappa^{1-d})\longrightarrow\g(V)_{A}
\]
by sending $a\otimes f$ to $\sum_{i\in A}a\otimes f_i(t)$ for any
$a\in V(d\,)$ and $f\in\Gamma(\Sigma\,;Q_{1},\,Q_{2},\ldots,Q_{N};\kappa^{1-d})$, where $\iota_{z_{i}}f(z_i)=\sum_{n\in\Z}c_{n}z_i^{n}$ is the
Laurent series expansion of the meromorphic function $f(z_{i})$ near $Q_{i}$ given by
$f=f(z_{i})(dz_{i})^{-\wt{a}+1}$, and $f_i(t)=\sum_{i\in\Z}c_{i}t^{i}\in\C((t))$. We denote the image of $j_{\bar{\Sigma}}$ by
$\g(V)_{\bar{\Sigma}}^{out}$; 
\[
\g(V)_{\bar{\Sigma}}^{out}=j_{\bar{\Sigma}}\left(\bigoplus_{d=0}^{\infty}
V(d\,)\otimes\Gamma(\Sigma\,;Q_{1},\,Q_{2},\ldots,Q_{N};\kappa^{1-d})\right).
\]
\begin{proposition}\label{Pro51} {\rm (\cite{Z2})} If $V=\bigoplus_{n=0}^\infty V(n),\, V(0)=\C\1$ 
is a quasi-primary generated vertex operator algebra 
and $\bar{\Sigma}$ is projective, then $\g(V)_{\bar{\Sigma}}^{out}$ is a Lie subalgebra of the Lie algebra $\g(V)_{A}$.
\end{proposition}
Let $\TS=(\Sigma\,;Q_{1},\,Q_{2},\ldots,Q_{N};z_{1},\,z_{2},\ldots,z_{N};W^{1},\, W^{2},\ldots,W^{N})$ be an $N$-labeled Riemann surface. 
We set  $W_{A}=\bigotimes_{i\in A}W^i$.  We denote by  $\rho_{Q_i}$ the action of $\g(V)_{i}$ on the $i$-th component of $W_{A}$, and set
$\rho_{\TS}=\bigoplus_{i\in A}\rho_{Q_i}$; the Lie algebra $\g(V)_{\Sigma}$ acts on $W_{A}$, and so the Lie subalgebra
$\g(V)_{\bar{\Sigma}}^{out}$. In other words we have a homomorphism $\rho_{\TS}:\g(V)_{\bar{\Sigma}}^{out}\rightarrow\End(W_A)$.

We now set  
\[
\GG(\TS)=\rho_{\TS}\left(j_{\bar{\Sigma}}\left(\bigoplus_{d=0}^{\infty}
\GG(V)(d\,)\otimes\Gamma(\Sigma\,;Q_{1},\,Q_{2},\ldots,Q_{N};\kappa^{1-d})\right)\right)\subset \End(W_A).
\] 
An element in $\GG(\TS)$ is called a \textit{quasi-global vertex operator}\,. If $V$ is quasi-primary generated and $\TS$ is projective, 
then $\GG(\TS)=\rho_{\TS}\left(\g(V)_{\bar{\Sigma}}^{out}\right)$, which is a Lie algebra by Proposition \ref{Pro51}. 
For any $a\otimes f$ we often denote $\rho_{\TS}(j_{\bar{\Sigma}}(a\otimes f))$ by $a(f,\TS)$ for simplicity.
The vector space $W_{A}$ is a module for the  Lie algebra $\GG(\TS)$.
We denote the space of coinvariants $W_{A}/\GG(\TS)W_{A}$ by $\qV(\TS)$. 

If $a$ is a primary vector, the quasi-global vertex operator $a(f,\TS)$ is called a \textit{global vertex operator} 
because it satisfies the transformation law like as $(1-|a|)$-differentials under coordinate changes.
Let $\PG({\TS})$ be the vector subspace of $\End(W_{A})$, which is linearly spanned by all global vertex operators 
and quasi-global vertex operators $\omega(f,\TS)$ for all $f\in\Gamma(\Sigma\,;Q_{1},\,Q_{2},\ldots,Q_{N};\kappa^{-1})$. 
 Then the {\it space of
covacua} is defined to be $\V(\TS)= W_{A}/\PG({\TS})W_{A}$. The main ingredient of this paper, the {\it space of vacua} or the {\it conformal
block} associated to the $N$-labeled Riemann surface $\TS$, is  defined to be $\V^{\dagger}(\TS)=\Hom_\C(W_{A}/\PG({\TS})W_{A},\C)$.

This definition of the space of covacua or the conformal is not convenient because in general $\PG({\TS})$ is not a Lie
subalgebra. However, due to the following theorem of Zhu, it suffices for us to consider the Lie algebra $\GG(\TS)$. 
\begin{theorem}\label{Thm1} {\rm (\cite[Theorem 5.2]{Z2})} 
Let $V=\bigoplus_{n=0}^\infty V(n)$ be a vertex operator algebra with $V(0)=\C\1$, and $\TS$ a
projective
$N$-labeled Riemann surface. Suppose that $V$ is a sum of highest weight 
modules for the Virasoro algebra. Then $\eta\in (W_A)^*$ belongs to
$\V^{\dagger}(\TS)$ if and only if $\eta(\GG(\TS)W_{A})=0$, i.e., there is a natural isomorphism as vector spaces
\[
\V^{\dagger}(\TS)\cong \qV(\TS)^*,
\]
where $\V^{\dagger}(\TS)$ is identified with the subset $\{\eta\in (W_{A})^{*}|\eta(\PG({\TS})W_{A})=0\}$ of $(W_{A})^{*}$.
\end{theorem}

In section \ref{sec:fin} a filtration on $\g(V)_{\bar{\Sigma}}^{out}$ being introduced here plays a very important role.
Let us start with an $N$-pointed Riemann surface 
$\bar{\Sigma}=(\Sigma\,;Q_{1},\,Q_{2},\ldots,Q_{N};z_{1},\,z_{2},\ldots,z_{N})$. 
For a given  meromorphic differential $f$ on $\Sigma$ whose poles are located at $Q_{1},\,Q_{2},\ldots,Q_{N}$ with order
$a_1,\,a_2,\ldots,a_N$, we define the order of $f$ by
\[
\divi{f}=\mathop{\max}\{a_1,\,a_2,\ldots,a_N\}.
\]
The filtration $\F_{p}\g(V)_{\bar{\Sigma}}^{out}\,(p\in\N)$ on $\g(V)_{\bar{\Sigma}}^{out}$ is defined by
\begin{equation}
\F_{p}\g(V)_{\bar{\Sigma}}^{out}=\har_{\C}\{j_{\bar{\Sigma}}(a\otimes f)\;|\;\wt{a}-1+\divi{f}\leq{p}\}.\label{Fil1}
\end{equation}
Then the Lie algebra $\g(V)_{\bar{\Sigma}}^{out}$ becomes a filtered Lie algebra.

We next define a filtration on $\g(V)_{\bar{\Sigma}}^{out}$-module $W_{A}$. 
We first recall that any $V$-module $W$ is a direct sum of $V$-modules of the form
$\bigoplus_{d=0}^{\infty} W(\lambda_{i}+d)$ $(i\in I)$ with lowest weight $\lambda_{i}$ 
such that $\lambda_i-\lambda_j\not\in\Z$ for some index set $I$.  We set $W_{d}=\bigoplus_{i\in I}W(\lambda_{i}+d)$ so that
$W=\bigoplus_{d=0}^\infty W_d$. The filtration $\F_pW_{A}\,(p\in\N)$ on $W_{A}$ is defined by
\[
{\F}_{p}W_{A}=\bigoplus_{0\leq d\leq p}W_{A,d},\quad W_{A,d}=\sum_{d_1+\cdots+d_N=d}W^{1}_{d_{1}}\otimes \cdots\otimes W^{N}_{d_{N}}.
\]
The $\g(V)_{\bar{\Sigma}}^{out}$-module $W_{A}$ becomes a filtered $\g(V)_{\bar{\Sigma}}^{out}$-module by this filtration.
 
Let $\F_{p}\qV(\TS)\,(p\in\N)$ be the induced filtration on $\qV(\TS)$, i.e, 
\[
\F_{p}\qV(\TS):=s(\F_{p}W_{A})=(\F_{p}W_{A}+\GG(\TS)W_{A})/\GG(\TS)W_{A},
\]
where $s$ is the natural projection $s:W_{A}\to\qV(\TS)$. We have the canonical surjection
\[
\pi:W_{A}=\bigoplus_{p=0}^\infty W_{A,p}\longrightarrow\gr_\bullet\qV(\TS):=\bigoplus_{p=0}^{\infty}\gr_{p}\qV( 
\TS),
\]
which is defined by $\pi(w)=s(w)+\F_{p-1}\V({\TS})\in\gr_{p}\V(\TS)$ for $w\in W_{A,p}$, 
where $\gr_{p}\,\qV(\TS):=\F_{p}\qV(\TS)/\F_{p-1}\qV(\TS)$.

\section{Finiteness of conformal blocks}\label{sec:fin}
We will prove the finiteness of conformal blocks over projective $N$-labeled Riemann surfaces $\widetilde{\Sigma}$\,: 
$\widetilde{\Sigma}=(\Sigma\,; Q_1,\,\ldots,Q_N\,;z_1,\ldots,z_N\,; W^1,\,\ldots, W^N)$. 
For a proof we basically follow the argument in \cite{AN}; however, we need to remedy the difficulties arising from the lack of
global meromorphic sections with lower order poles of $\kappa^{1-n}$ for positive integer $n$, i.e., Weierstrass gaps, 
which do not appear in the case of the projective line.

In \cite{AN} we explicitly constructed global meromorphic sections of $\kappa^{1-n}\,(n\geq1)$ for the canonical bundle $\kappa$ of the projective
line,  which have poles of desired orders at the prescribed point, and are holomorphic elsewhere. 
Using Riemann-Roch theorem for any compact Riemann surface $\Sigma$,
we see that there exists a global meromorphic section which has the pole at $Q\in\Sigma$, and is holomorphic on
$\Sigma\setminus\{Q\}$; however in general the order is large 
so that we  are not able to find such a meromorphic section having the lower order poles at $Q$. 

\begin{lemma}\label{Thm4}
Let $\Sigma$ be a  compact Riemann surface of genus $g$. We fix a point $Q\in\Sigma$ and a positive integer $n\in\Z_{>0}$. 
\vskip 1ex
\noindent
{\rm (1)} There exists a nontrivial global meromorphic section $f$ of $\kappa^{1-n}$, which has a pole at $Q$ and is holomorphic on 
$\Sigma\setminus\{Q\}$.

\vskip 1ex
\noindent
{\rm (2)} Let $\nu$ be the order of the pole at $Q$ of the global meromorphic section $f$ in {\rm (1)}.  
Set $M=\nu+2g$. Then for any $m\geq M$, there exists a global meromorphic section of $\kappa^{1-n}$,
which has a pole of order $m$ at $Q$ and is holomorphic on $\Sigma\setminus\{Q\}$.
\end{lemma}
\begin{proof}
The assertion (1) is found in \cite[Theorem 29.16, page 225]{F} or it is directly proved by using Riemann-Roch theorem. 
By Weierstrass gap theorem (\cite[page 202]{B}), for any $i\in\N$ there exists a meromorphic function $h$ on $\Sigma$ 
such that $h$ has a pole of order $i+2g$ at $Q$ and holomorphic of $\Sigma\setminus\{Q\}$. Then $hf$ has the pole at $Q$ of order $2g+\nu+i$.
\end{proof}

Let $U$ be a subspace of $V$, which is linearly spanned by finitely many homogeneous quasi-primary vectors. Let $r_U$ be the maximal number
of weight of homogeneous vectors of $U$. Using Lemma \ref{Thm4} we can find a positive integer $M_U$ such that for any $n\leq r_U,\,m\geq M_U$ and
$i\in A$, there exists a global meromorphic section over $\kappa^{1-n}$ which has a pole of order $m$ at $Q_{i}$ and is holomorphic of
$\Sigma\setminus\{Q_{i}\}$. 

We denote by  $C_m(U,W)\,(m\geq2)$ the subspace of $W$, which is linearly spanned  by $a(-n)w$ for all $a\in U,\,w\in
W$ and $n\geq m$.

\begin{lemma}\label{Lemg1}
Let $V=\bigoplus_{n=0}^\infty V(n),\, V(0)=\C\1$ be a quasi-primary generated vertex operator algebra.
Let $U$ be a subspace of $V$, which is linearly spanned by finitely many quasi-primary vectors. 
Then for any $i\in A$ the set $W^{1}\otimes\cdots\otimes C_{M_U}(U, W^{i})\otimes\cdots\otimes W^{N}$ 
is in the kernel of the surjective linear map
$\pi:W_{A}\to\gr_\bullet\,\qV(\TS)$. In particular, the map $\pi$ induces a surjective linear map
\[
\pi:W^{1}/C_{M_U}(U,W^{1})\otimes\cdots\otimes W^{N}/C_{M_U}(U,W^{N})\longrightarrow\gr_\bullet\,\qV(\TS).
\]
\end{lemma}
\begin{proof}
It suffices to show that for a homogeneous $a\in U$
\[
\pi(w_{1}\otimes\cdots\otimes a(-m)w_{i}\otimes\cdots\otimes w_{N})=0 
\mbox{ for all $w_{i}\in W^{i}_{d_{i}}$ and $m\geq M_U$.}
\]
By definition of $M_U$, for any $m\geq M_{U}$, there exists 
$f\in\Gamma(\Sigma\,;Q_{i}\,;\kappa^{1-\wt{a}})$ which has the Laurent 
series expansion
$\iota_{z_{i}}f=z_{i}^{-m}+\sum_{\ell>-m}c_{\ell}z_{i}^{\ell}$ at $Q_i$\,. 
Then by definition of quasi-global vertex operators we see that
\begin{align*}
&a(f,\TS)(w_{1}\otimes\cdots\otimes w_{N})\\
&\qquad=w_{1}\otimes\cdots\otimes a(-m)w_{i}\otimes\cdots\otimes 
w_{N}+\sum_{\ell>-m}c_{\ell}w_{1}\otimes\cdots\otimes
a(\ell)w_{i}\otimes\cdots\otimes w_{N}\\
&\qquad\qquad+\sum_{\substack{j\in A\\j\neq i}}w_{1}\otimes\cdots\otimes
(\Res_{z_{j}}Y(a,z_{j})\iota_{z_{j}}f)w_{j}\otimes\cdots\otimes w_{N}\\
&\qquad=w_{1}\otimes\cdots\otimes a(-m)w_{i}\otimes\cdots\otimes 
w_{N}+\mbox{lower weight (degree) terms}\in\GG(\TS)W_{A},
\end{align*}
where we have used the fact that $f$ is holomorphic at $Q_j\,(j\neq i)$. This implies that
\[
w_{1}\otimes\cdots\otimes a(-m)w_{i}\otimes\cdots\otimes w_{N}\in\F_{\sum 
d_{i}+\wt{a}+m-2}W_{A}+\GG(\TS)W_{A}.
\]
Since $w_{1}\otimes\cdots\otimes a(-m)w_{i}\otimes\cdots\otimes w_{N}\in W_{A,\sum\!d_{i}+\wt{a}+m-1}$,
this element belongs to the kernel of the map $\pi$.
\end{proof}

Let $U$ be a graded subspace of $V$ such that $V=U\oplus C_2(V)$. 
Recall that $V$ is linearly spanned by vectors $\alpha_{1}(-n_{1})\cdots\alpha_{r}(-n_{r})\1\,(\alpha_i\in U)$ with 
$n_{1}>\ldots>n_{r}>0$. For any positive integers $m$ and $q$, we set
\[
C_{m,q}(W)=\{\,(\alpha_{1}(-n_{1})\cdots\alpha_{r}(-n_{r})\1)(-p)w\,|\,m\geq n_{1}>\ldots>n_{r}>0,\,\alpha_{i}\in U,\,p\geq q\,\}.
\]

\begin{lemma}\label{Lemg2}
Let $m,\,q$ be positive integers. Then $C_q(W)\subset C_{m,q}(W)+C_{m}(U,W)$.
\end{lemma}
\begin{proof} By Proposition \ref{Pro7} it suffices to show that for any $\alpha_{i}\in U$ and $n_{1}>\ldots>n_{r}>0$
\begin{equation}\label{equ1}
(\alpha_{1}(-n_{1})\cdots\alpha_{r}(-n_{r})\1)(-p)w\in C_{m,q}(W)+C_{m}(U,W)
\end{equation}
for all $p\geq q$; we can assume that $n_1>m$ by definition of $C_{m,q}(W)$. 

We now see that 
\[
(\alpha(-n)\1)(-p)w=(-1)^{p-1}\binom{-n}{p-1}\alpha(-n-p+1)w\in C_{m}(U,W)\mbox{ for all $\alpha\in U$},
\]
which proves the case $r=1$.
Suppose that (\ref{equ1}) holds for any $r<r_{0}$ for some $r_0\geq2$. 
We set $\beta=\alpha_{2}(-n_{2})\cdots\alpha_{r_0}(-n_{r_{0}})\1$, and use the  associativity formula (\ref{asso1}) to get
\begin{multline*}
(\alpha_{1}(-n_{1})\beta)(-p)w\\
=\sum_{i=0}^{\infty}\binom{-n_1}{i}(-1)^{i}\left(\alpha_{1}(-n_{1}-i)\beta(-p+i)w-( 
-1)^{n_{1}}\beta(-n_{1}-p-i)\alpha_{1}(i)w\right).
\end{multline*}
Then the first term and the second of the right hand side  belongs to $C_m(U,V)$ and $C_{m,q}(W)+C_{m}(U,W)$, respectively, 
where we use inductive hypothesis to the second term.
\end{proof}

\begin{lemma}\label{Lemg3}
Let $m$ be a positive integer. For any  $a\in V\,(|a|\geq1)$ and $w\in W$ 
we have $a(-q)w\in C_{m}(U,W)$ for all $q$ such that $q\geq m\wt{a}$.
\end{lemma}
\begin{proof} 
If $\wt{a}=1$ then $a\in U$ 
because $V(1)\cap C_{2}(V)=\{0\}$. Hence $a(-q)w\in C_{m}(U,W)$
for any $q\geq m(=\wt{a}m)$. Suppose that $\wt{a}>1$. If $a\in U$ then 
$a(-q)w\in C_{m}(U,W)$ for any $q\geq \wt{a}m(>\!m)$.
We can now assume that $a\in C_{2}(V)$.
Suppose that $(0\neq)a=b'(-2)c$ for some $b'$ and $c$.
Then $a=(L_{-1}b')(-1)c=b(-1)c$ and $1\leq\wt{b}\leq \wt{a},\,\wt{c}<\wt{a}$.

By the associativity formula we have
\[
a(-q)w=(b(-1)c)(-q)w=\sum_{i=0}^{\infty}\big(b(-1-i)c(-q+i)w+c(-1-q-i)b(i)w
\big)
\]
for any $q\in\Z$.  Since $q\geq\wt{a}m>\wt{c}m$ using inductive hypothesis to the element $c$ we see that $c(-1-q-i)b(i)w\in C_{m}(U,W)$ 
for $i\geq 0$. We will show that $b(-1-i)c(-q+i)w\in C_{m}(U,W)$ for any $q\geq\wt{a}m$ and $i\geq 0$. 
If $i\geq\wt{b}m$ then  $b(-1-i)c(-q+i)w\in C_{m}(U,W)$ for any $q\in\Z$ by inductive hypothesis. 
Otherwise, i.e., $i<\wt{b}m$,
recall the commutator formula
\[
b(-1-i)c(-q+i)w=c(-q+i)b(-1-i)w+\sum_{j=0}^{\infty}\binom{-1-i}{j}(b(j)c)(-1 
-q-j)w.
\]
Now since $q-i\geq \wt{a}m-\wt{b}m+1=(\wt{a}-\wt{b})m+1=\wt{c}m+1>\wt{c}m$, using inductive hypothesis we see that 
the first term $c(-q+i)b(-1-i)w$ belongs to $C_{m}(U,W)$. Finally since $\wt{a}>\wt{b(j)c}$ for any $j\geq0$, 
inductive hypothesis shows that $(b(j)c)(-1-q-j)w\in C_{m}(U,W)$.
\end{proof}

\begin{proposition}\label{Prog4}
Let $V=\bigoplus_{n=0}^\infty V(n),\,V(0)=\C\1$ be a $C_{2}$-finite vertex operator algebra, 
and $U$ be any finite dimensional graded subspace of $V$ such that $V=U\oplus C_2(V)$. Let $m$ be a positive integer. 
Then there exists a positive integer $k$ such that $C_{k}(W)\subset C_{m}(U,W)$.
\end{proposition}
\begin{proof}
By Lemma \ref{Lemg2} we know that $C_k(W)\subset C_{m,k}(W)+C_{m}(U,W)$ for any positive integer $k$. 
Then it suffices to show that there exists a positive integer $k$ such that $C_{m,k}(W)\subset C_{m}(U,W)$.
Let $s_U$ be the maximal number of weights of homogeneous elements in $U$. 
For any  $\alpha_{i}\in U\,(1\leq i\leq r)$ and $m\geq n_{1}>\ldots>n_{r}>0$ we see that 
\begin{align*}
\wt{\alpha_{1}(-n_{1})\cdots\alpha_{r}(-n_{r})\1}&=\sum_{i=1}^{r}(\wt{\alpha_{i}}-1)+\sum_{i=1}^{r}n_{i}\\
&\leq r(s_{U}-1)+\sum_{i=1}^{r}(m-i+1)=-\frac{1}{2}r^{2}+r\left(s_{U}+m-\frac{1}{2}\right),
\end{align*}
which is bounded from above by some positive integer $k_0\geq\left(s_{U}+m-1/2\right)^2/2$. 
We note that the constant $k_0$ depends only on $m$ and $U$. Setting $k=k_0m$ we get $C_{m,k}(W)\subset C_{m}(U,W)$ 
by Lemma \ref{Lemg3} because any element in $C_{m,k}(W)$ is a linear combination of $a(-p)w$ 
for $\wt{a}\leq k_{0},\, w\in W$ and $p\geq k\,(\geq m|a|)$.

\end{proof}

\begin{proposition}\label{Prog5}
Let $V=\bigoplus_{n=0}^\infty V(n),\,V(0)=\C\1$ be a $C_{2}$-finite vertex operator algebra and $W$ a weak $V$-module. 
If the module $W$ is $B_{1}$-finite then $W/C_{m}(U,W)$ is of finite dimensional for any $m>0$.
\end{proposition}
\begin{proof}
By Proposition \ref{Prog4} there exists a positive integer $k$ such that $C_{k}(W)\subset C_{m}(U,W)$. Since $V$ is $C_{2}$-finite and $W$ is
$B_{1}$-finite, $W$ is $C_{k}$-finite by Theorem \ref{Thm3}. Thus $W/C_{k}(W)$ is of finite dimensional, and so is $W/C_{m}(U,W)$.
\end{proof}

\begin{theorem}\label{Thm5}
Let $V=\bigoplus_{n=0}^\infty V(n),\,V(0)=\C\1$ be a quasi-primary generated, $C_{2}$-finite vertex operator algebra 
such that $V$ is a sum of highest weight modules for the Virasoro algebra,
and let $\TS=(\Sigma\,;Q_{1},\,\ldots,Q_{N};z_{1},\ldots,z_{N};W^{1},\ldots,W^{N})$ be  a projective $N$-labeled Riemann surface. 
If all $V$-modules $W^{i}\,(i\in A)$ are $B_{1}$-finite, then the conformal block $\V^{\dagger}(\TS)$ is of finite dimensional.
\end{theorem}
\begin{proof} Let $U$ be a finite dimensional graded subspace of $V$ such that $V=U\oplus C_{2}(V)$. 
Since $V$ is quasi-primary generated
any vectors from $U$ are linear combinations of $L_{-1}^ia$ for some $i\in\N$ and quasi-primary vectors $a$.
Then we can further assume that any elements of $U$ are quasi-primary because $L_{-1}a=a(-2)\1\in C_2(V)$ for any $a\in V$.

By Theorem \ref{Thm1} it suffices to prove that $\qV(\TS)$ is of finite dimensional. We set $M=M_U>0$. The constant $M_U$ is defined in the
paragraph just before Lemma
\ref{Lemg1}; recall that in order to define $M_U$ we assume that $U$ is linearly spanned by quasi-primary vectors.  By Proposition \ref{Prog5}
$W^{i}/C_{M}(U,W^{i})$ is of finite dimensional for any $i\in A$.  Thus Lemma \ref{Lemg1} shows that
$\gr_\bullet\,\qV(\TS)$ is of finite dimensional and so is $\qV(\TS)$. 

\end{proof}


\section{Examples}
We present several examples of $C_{2}$-finite vertex operator algebras; affine vertex operator algebras (with positive integral level $k$), 
Virasoro vertex operator algebras (with minimal central charge $c_{p,q}$) and lattice vertex operator algebras. 
We will prove that all irreducible modules for these vertex operator algebras are $B_1$-finite. 
Then we see that those modules are all $C_2$-finite by Theorem \ref{Thm3}.
In fact $C_2$-finiteness for irreducible modules for affine and Virasoro vertex operator algebra is well known (cf. \cite{DLM} for affine case, and
\cite{FF2} for Virasoro case). The $C_2$-finiteness for irreducible modules for lattice vertex operator algebras seem to be known, 
though we are not able to find any published materials so far.

In order to prove the $B_1$-finiteness in Virasoro case we follow the same argument being used in the proof of
$C_2$-finiteness, however, we will see that verifying $B_1$-finiteness is much easier than $C_2$-finiteness.

\begin{example}[Affine vertex operator algebras]
Let $\hat{\g}=\C[t,t^{-1}]\otimes\g\oplus\C c\oplus\C d$ be an affine Lie algebra where $\g$ is a finite dimensional simple Lie algebra.
We denote by $\{\Lambda_{0},\ldots,\Lambda_{n}\}$ the set of fundamental weights for $\hat{\g}$, and by $P^{k}_{+}$ the set of all level $k$
dominant integral weights. We denote the irreducible highest weight module of $\hat{\g}$ with highest weight $\Lambda$ by $L(\Lambda)$. 
It is known that if $k\neq -h^\vee,0$ where $h^\vee$ is the dual Coxeter number of $\hat{\g}$, then $L_k=L(k\Lambda_{0})$ is a vertex operator
algebra. Moreover, if $k$ is positive integer any irreducible $L_k$-module is realized as an irreducible $\hat{\g}$-module
$L(\Lambda)$ for some $\Lambda\in P^{k}_{+}$ (see \cite{FZ}). The $C_{2}$-finiteness of $L_{k}$ is known (\cite{Z1}, \cite{DLM}).

We now prove the $B_{1}$-finiteness of irreducible $L_{k}$-modules. Since $L(\Lambda)$ is linearly spanned by vectors 
$a_{1}(-n_{1})\cdots a_{r}(-n_{r})v$ with $n_{i}>0$, $a_{i}\in L_{k}(1)(\cong\g)$ and $v\in V_{\Lambda}$, 
where $V_\Lambda$ is the irreducible  highest weight module for $\mathfrak{g}$ with the highest weight $\bar{\Lambda}$ 
and highest weight vector $v_{\bar{\Lambda}}$, where $\bar{\Lambda}$ is the classical part of $\Lambda$. 
We now see that $L(\Lambda)=V_{\Lambda}+B_{1}(L(\Lambda))$, and that $L(\Lambda)$ is $B_{1}$-finite 
because $V_{\Lambda}$ is of finite dimensional.
\end{example}
\begin{example}[Lattice vertex operator algebras]
Let $L$ be an even lattice of finite rank with a nondegenerate positive definite symmetric $\Z$-bilinear form
$\langle\,\cdot\,|\,\cdot\,\rangle$. We set $\h=\C\otimes_{\Z}L$ and
$\hat{\h}=\C[t,t^{-1}]\otimes_{\C}\h\oplus\C K$; the latter is the affinization of the abelian Lie algebra $\h$. 
Let $L^{\circ}$ be the dual lattice of $L$, and $\C[L^{\circ}]=\bigoplus_{\beta\in L^{\circ}}\C\,e_{\beta}$ be the twisted group algebra of
$L^{\circ}$ with some cocycle which represents the central extension of $L^\circ$. For any subset $M$ of $L^{\circ}$ we write
$\C[M]=\bigoplus_{\beta\in M}\C\,e_{\beta}$,  and set $V_{M}=U(\hat{\mathfrak{h}}^-)\otimes\C[M]$ where
$\hat{\mathfrak{h}}^-=t^{-1}\C[t^{-1}]\otimes_{\C}\h$ is a Lie subalgebra  of
$\hat{\mathfrak{h}}$. 
We note that the Lie algebra $\hat{\mathfrak{h}}$ naturally acts on $V_M$.
Then it is known that $V_{L}$ is a vertex operator algebra, and that the $V_{\lambda+L}$ for $\lambda\in L^{\circ}$ give all irreducible 
$V_{L}$-modules (see \cite{FLM}). The vertex operator associated to
$e_{\alpha}\,(\alpha\in L)$ is 
\[
Y(e_{\alpha},x)=\exp\left(\sum_{n=1}^{\infty}{\frac{\alpha(-n)}{n}}x^{n}\right)
\exp\left(-\sum_{n=1}^{\infty}{\frac{\alpha(n)}{n}}x^{-n}\right)e_{\alpha}x^ {\alpha(0)},\quad \alpha(n)=t^{n}\otimes\alpha
\]
where $e_{\alpha}$ acts on $\C[L^{\circ}]$  by the left
multiplication, and the action of $x^{\alpha(0)}$  on $V_{L^{\circ}}$ is defined by 
$x^{\alpha(0)}(u\otimes e_\mu)=x^{\langle\alpha|\mu\rangle}(u\otimes e_\mu)$ for all $\mu\in L^{\circ}$ and $u\in U(\hat{\mathfrak{h}}^-)$. 

We now prove that for any $\lambda\in L^{\circ}$ the irreducible $V_{L}$-module $V_{\lambda+L}$ is $B_{1}$-finite. 
We set
$\Gamma_{\lambda}=\{\,\beta\in L\,|\,\langle\beta-\alpha|\alpha+\lambda\rangle<0 \mbox{ for any $\alpha\in
L$  such that } \alpha\neq\beta,-\lambda\,\}$. The following lemma is due to  H.~Shimakura.
\begin{lemma}\label{lattice2}
Let $\lambda\in L^{\circ}$. Then $\Gamma_{\lambda}$ is a finite set.
\end{lemma}
\begin{proof} 
Let $\phi$ be the translation map on $L^{\circ}$  defined by $\phi(\gamma)=\gamma+\lambda$\,.  
We see that $\Gamma_{\lambda}=\{\,\beta\in L\,|\,\langle\delta|\beta+\lambda-\delta\rangle<0 \hbox{ for any $\delta\in L$ 
such that $\delta\neq0,\lambda+\beta$}\,\}$, and that 
\[
\phi(\Gamma_\lambda)=\{\,\gamma\in\lambda+L\,|\,\langle\delta|\gamma-\delta\rangle<0 \hbox{ for any $\delta\in L$ such that
$\delta\neq0,\gamma$}\}.
\]

Let $\alpha_{1},\ldots,\alpha_{\ell}$\, be a basis of $L$, and let $\Lambda_{1},\ldots,\Lambda_{\ell}$  be the basis of $L^\circ$ such
that $\langle\Lambda_i|\alpha_j\rangle=\delta_{ij}$. 
Let $\gamma\in\phi(\Gamma_\lambda)$ and $\gamma\neq\pm\alpha_i\,(1\leq i \leq\ell)$. By definition of $\phi(\Gamma_\lambda)$ we have
$\langle\gamma-(\pm\alpha_{i})|\pm\alpha_{i}\rangle<0\;(i=1,\ldots,\ell)$, and we see that
$\gamma\in\{\,\sum_{i=1}^{\ell}m_{i}\Lambda_{i}\,|\,m_{i}\in\Z\hbox{ and }|m_{i}|\leq\langle\alpha_{i}|\alpha_{i}\rangle\hbox{ for any $i$}\,\}$. 
Therefore, $\phi(\Gamma_\lambda)$ is a finite set, and so is
$\Gamma_{\lambda}$. 
\end{proof}

By definition we see that $V_{\lambda+L}=\sum_{\alpha\in L}\C e_{\lambda+\alpha}+B_{1}(V_{\lambda+L})$, and
for any $\alpha,\beta\in L$ we have
$e_{\beta-\alpha}(-\langle\beta-\alpha|\lambda+\alpha\rangle-1)e_{\lambda+\alpha}=\pm e_{\lambda+\beta}$. Hence we find
that $V_{\lambda+L}=\sum_{\alpha\in \Gamma_{\lambda}}\C e_{\lambda+\alpha}+B_{1}(V_{\lambda+L})$. Then Lemma \ref{lattice2} shows that
$V_{\lambda+L}$ is $B_{1}$-finite.
\end{example}

\begin{example}[The Virasoro vertex operator algebras] 
Let $\Vir=\bigoplus_{n\in\Z}\C L_{n}\oplus\C C$ be the Virasoro algebra. 
Let $M(c,h)$ be the Verma module for the Virasoro algebra with a highest weight $h\in\C$ and central charge $c\in\C$.
We denote by $v_{h,c}$ the highest weight vector, i.e.,
$v_{h,c}$ satisfies  $L_{n}v_{h,c}=\delta_{n,0}hv_{h,c}\,(n\geq0)$ and
$Cv_{h,c}=cv_{h,c}$. The Verma module $M(c,h)$ is a rank one free $U(\Vir^-)$-module with the generator $v_{h,c}$ where 
$\Vir^-=\bigoplus_{n\in\Z_{>0}}\C L_{-n}$. Let
$L(c,h)$ be the irreducible quotient of $M(c,h)$. Then it is known that $L(c,0)$ is a vertex operator algebra.

Let $p,q$ be coprime positive integers.  We set $c_{p,q}=1-6(p-q)^{2}/pq$. For any integers $r$ and $s$ 
such that $1\leq r<q,\,1\leq s<p$ we denote $h_{p,q;r,s}=((rp-sq)^{2}-(p-q)^{2})/4pq$. 
Then any irreducible $L(c_{p,q},0)$-module is isomorphic to $L(c_{p,q},h_{p,q;r,s})$.  We prove: 
\begin{proposition}\label{Pro71}
$L(c_{p,q},0)$ is $C_{2}$-finite and any irreducible $L(c_{p,q},0)$-module $L(c_{p,q},h_{p,q;r,s})$ is $B_{1}$-finite.
\end{proposition}

In order to prove the proposition we recall several properties of singular vectors $v$ in $M(c_{p,q},h_{p,q;r,s})$, i.e., 
$L_{n}v=v$ for $n>0$. 
It is known that for any positive integers $r$ and $s$ satisfying $1\leq r<q$ and $1\leq s<p$ there exists a unique singular
vector $u_{r,s}\in M(c_{p,q},h_{p,q;r,s})$ such that $L_0u_{r,s}=(h_{p,q;r,s}+rs)u_{r,s}$  up to scalar multiples.  The explicit form of the
singular vector is not known, but we have the partial formula which expresses this singular vector as explained below.

Let us fix central charge $c=c_{p,q}$ and highest weight $h=h_{p,q;r,s}$.
Let $\Vir^{\leq-3}=\bigoplus_{n\geq 3}\C L_{-n}$ which is a Lie subalgebra of $\Vir$.
Let $\phi:\C[x,y]\to M(c,h)/\Vir^{\leq-3}M(c,h)$ be the linear isomorphism defined by 
\[
x^{i}y^{j}\longmapsto
L_{-2}^{j}L_{-1}^{i}v_{h_{r,s},c}+\Vir^{\leq-3}M(c,h).
\]
Let $f:M(c,h)\to M(c,h)/\Vir^{\leq-3}M(c,h)$ be the canonical projection.
We define $\pi=\phi^{-1}\circ f$. The following proposition is proved in \cite{FF1}:
\begin{proposition}\label{Pro8} 
Let $c=c_{p,q}$ and $h=h_{p,q;r,s}$. Let $u_{r,s}\in M(c,h)$ be the singular vector such that $L_0u_{r,s}=(h+rs)u_{r,s}$.
Then $\pi(u_{r,s})=\alpha F_{r,s}(x,y;p/q)$ for some nonzero constant $\alpha$, where $F_{r,s}(x,y;t)$ is a
polynomial of $\C[x,y,t,t^{-1}]$ given by 
\[ 
F_{r,s}(x,y;t)^{2}=\prod_{k=0}^{r-1}\prod_{\ell=0}^{s-1}(x^{2}-\{(r-2k-1)t^{1/2}-(s-2\ell-1)t^{-1/2}\}^{2}y).
\]
\end{proposition} 

We now can prove Proposition \ref{Pro71}.
\begin{proof}[Proof of Proposition \ref{Pro7}.] 
The canonical projection $\psi:M(c,h)\to L(c,h)$ maps the subspace $\Vir^{\leq-3}M(c,h)$ into $C_{2}(L(c,h))$. 
Any singular vectors in $M(c,h)$ are in the kernel of this map. We note that $h_{p,q;r,s}=h_{p,q;q-r,p-s}$, in particular there exists
a singular vector $u_{q-r,p-s}$ such that $L_0u_{q-r,p-s}=(h+(p-s)(q-r))u_{q-r,p-s}$.  We see that
the composition of $\phi$ and $\psi$ induces a surjective linear map  
\begin{align}\label{equ2}
\psi:\C[x,y]/(F_{r,s}(x,y;p/q),F_{q-r,p-s}(x,y;p/q))\longrightarrow L(c,h)/C_{2}(L(c,h)).
\end{align} 
First we find that $F_{1,1}(x,y;p/q)=x$ and $F_{q-1,p-1}(x,y;p/q)\equiv \alpha' y^{(p-1)(q-1)/2} \mod (x)$ for some nonzero $\alpha'\in\C$. 
Thus $\C[x,y]/(F_{1,1}(x,y;p/q),F_{q-1,p-1}(x,y;p/q))$ is of finite dimensional, which show that  $L(c,0)$ is $C_{2}$-finite. 

Finally we prove that any irreducible module $L(c,h)$ is $B_{1}$-finite. We find that the surjective map (\ref{equ2}) induces a surjective
map  
\begin{align}
\C[x,y]/(y, F_{r,s}(x,y;p/q),F_{q-r,p-s}(x,y;p/q))\to L(c,h)/B_{1}(L(c,h)).
\end{align} 
Since $\C[x,y]/(y, F_{r,s}(x,y;p/q),F_{q-r,p-s}(x,y;p/q))\cong\C[x]/(x^{rs},x^{(q-r)(p-s)})$ is of finite dimensional we see that $L(c,h)$ is
$B_{1}$-finite for any $1\leq r<q$ and $1\leq s<p$.  
\end{proof}
\end{example}

\end{document}